# Cutting Cakes Correctly

by Theodore P. Hill
School of Mathematics, Georgia Institute of Technology, Atlanta, GA 30332-0160
hill@math.gatech.edu

The article "Better Ways to Cut a Cake" [BJK] by Brams, Jones and Klamler in the December 2006 issue of these *Notices* is receiving widespread media attention in *Scientific American Science News*, *Science Daily*, and *the Discovery Channel*, among others, as well as in our own Society's promotional site *AMS in the News*. Unfortunately, [BJK] contains serious mathematical errors, some of which will be summarized here. Further details may be found in [TPH].

In the first section of [BJK], the authors state (p 1314) that in problems of fair division of a divisible good, "the well-known 2-person, 1-cut cake-cutting procedure 'I cut, you choose'" is Pareto-optimal, that is, "There is no other allocation that is better for one person and at least as good for the other." Cut-and-choose is not even Pareto optimal among (*n*-1)-cut procedures, a weaker form of Pareto optimality called "C-efficient" [BT p 149-150], as the following simple example shows.

**Counterexample 1**. The "cake" is the unit square, and player 1 values only the top half of the cake and player 2 only the bottom half (and on those portions, the values are uniformly distributed). If player 1 is the cutter, and cuts vertically, his uniquely optimal cut-and-choose solution is to bisect the cake exactly, in which each player receives a portion he values exactly Ω. Or if player 1 cuts horizontally, his uniquely optimal risk-adverse cut-and-choose point is the line y = æ, in which case he receives a portion he values at Ω the cake, and player 2 chooses the bottom portion and receives a portion he values at 100% of the cake. But an allocation of the top half of the cake to player 1, and the bottom half to player 2 is at least as good for player 2 in both cases, and is strictly better for player 1, so cut-and-choose is not Pareto optimal in either direction.

If the cake is the unit interval and the value measures are all continuous as well as mutually absolutely continuous, then cut-and-choose is Pareto optimal. That conclusion may fail if the measures are not continuous, as is easily seen by looking at the case where all players place all the value of the cake on the same single point (and hence are mutually absolutely continuous). And if all the measures are continuous but not mutually absolutely continuous (which was not assumed anywhere in [BT]), the statement quoted from [BT] in [BJK, footnote 3 p 1318] "an envy-free allocation that uses *n* -1 parallel cuts is always efficient [i.e., Pareto optimal]", as well as the corresponding Proposition 7.1 of [BT, p 150], are not true.

**Counterexample 2**. The cake is the unit interval; player 1 values it uniformly, and player 2 values only the left- and right-most quarters of the interval, and values them equally and uniformly. (In other words, the probability density function (pdf) representing player 1's value is a.s. constant 1 on [0,1], and that of player 2 is a.s. constant 2 on [0, º] and on [3/4,1], and zero otherwise.) If player 1 is the cutter, his unique cut point is at x = Ω, and each player will receive a portion he values at exactly Ω. The allocation of the interval [0, º] to player 2 and the rest to player 1, however, gives player 1 a portion he values æ, and player 2 a portion he values Ω again, so cut-and-choose (which is an envy-free allocation for 2 players) is not Pareto optimal.

The new cake-cutting procedure described in [BJK], Surplus Procedure (SP), is not well defined. The integrals defining SP [BJK p 1315] need not exist, if, for example, the measures do not have densities. Also, if a player's value measure does not have a unique median, the defining cut-point in SP is not unique. Similarly, without additional assumptions, the new Equitability Procedure (EP) in [BJK] is neither well defined nor constructive, since the underlying systems of $n$-1 integral equations in $n$-1 unknowns may not have solutions, and those that do may not have solutions in closed form (cf. [TPH]).

The authors claim that both EP and SP are Pareto optimal [BJK, pp 1318, 1320], but neither procedure is Pareto optimal as defined in [BJK, (2) p 1314] and as defined in [BT, page 44]; for counterexamples, see [TPH]. The underlying reason is that both EP and SP allocate contiguous portions to each player, and as noted in [BT, p 149], "satisfying contiguity may be inconsistent with satisfying efficiency [Pareto optimality]".

The authors [BJK, 1316] define an allocation procedure to be *strategy-vulnerable* if a "player can, by misrepresenting its value function, assuredly do better, whatever the value function of the other player"; and otherwise the procedure is called *strategy-proof*. The second conclusion of [BJK,Theorem 1] states "any procedure that makes $e$ the cut-point is strategy vulnerable". That conclusion is false, as the next example shows.

**Counterexample 3.** Suppose player 1 has true value measure *v*. Every procedure allocates disjoint subsets of the cake, one to each player, and if the players' values happen to be identical, at least one player receives a portion he values no more than *1/n*. That player could be player 1, so he will not do "assuredly better" than his fair share of *1/n*. Hence every fair allocation procedure is strategy-proof, including the one that makes *e* the cut-point, contradicting the conclusion of Theorem 1 (and showing that the first conclusion of [BJK Theorem 1] and [BJK, Theorem 2] are trivial).

The fifth sentence in the proof of [BJK, Theorem 3] says "By moving all players' marks rightward … one can give each player an equal amount greater than 1/n". This is fallacious (for counterexample see [TPH]). The rest of the proof of [BJK, Theorem 3] is incomplete, since it only proves a claim about the moving-knife procedure, whereas the desired conclusion concerns EP.

**Acknowledgement.** The author is grateful to Professor Kent E. Morrison and two anonymous colleagues for valuable comments and suggestions.